\documentclass[12pt,a4]{article}%
\usepackage{amsmath}
\usepackage{graphicx}
\usepackage{hyperref}
\usepackage{eurosym}
\usepackage{amsfonts}
\usepackage{amssymb}%
\providecommand{\U}[1]{\protect\rule{.1in}{.1in}}
\providecommand{\U}[1]{\protect\rule{.1in}{.1in}}

\begin{document}

\begin{center}
{\Large Stochastic Fixed Points\ and Nonlinear}\medskip

{\Large Perron-Frobenius Theorem}\footnote{The authors are grateful to B.M.
Gurevich, V.I. Oseledets, K.R. Schenk-Hopp\'{e} and participants of the
seminar of the Dobrushin Laboratory at the Institute for Information
Transmission Problems in Moscow for helpful comments and fruitful
discussions.}\bigskip

E. Babaei\footnote{Economics, University of Manchester, Oxford Road,
Manchester M13 9PL, UK. E-mail:
esmaeil.babaeikhezerloo@postgrad.manchester.ac.uk.}, I.V.
Evstigneev\footnote{Economics, University of Manchester, Oxford Road,
Manchester M13 9PL, UK. E-mail: igor.evstigneev@manchester.ac.uk.}, and S.A.
Pirogov\footnote{Moscow University and IITP, Academy of Sciences of Russia,
Bolshoy Karetny 19-1, Moscow, 127051, Russia. E-mail: pirogov@bk.ru .}\bigskip
\end{center}

\begin{quotation}
{\small \textbf{Abstract: }We provide conditions for the existence of
measurable solutions to the equation }$\xi(T\omega)=f(\omega,\xi(\omega
))${\small , where }$T:\Omega\rightarrow\Omega$ {\small is an automorphism of
the probability space }$\Omega$ {\small and }$f(\omega,\cdot)$ {\small is a
strictly non-expansive mapping. We use results of this kind to establish a
stochastic nonlinear analogue of the Perron-Frobenius theorem on eigenvalues
and eigenvectors of a positive matrix. We consider a random mapping }%
$D(\omega)$ {\small of a random closed cone }$K(\omega)$ {\small in a
finite-dimensional linear space into the cone }$K(T\omega)${\small . Under
assumptions of monotonicity and homogeneity of }$D(\omega)${\small , we prove
the existence of scalar and vector measurable functions }$\alpha(\omega)>0$
{\small and }$x(\omega)\in K(\omega)$ {\small satisfying the equation }%
$\alpha(\omega)x(T\omega)=D(\omega)x(\omega)$ {\small almost surely.}\bigskip
\end{quotation}

{\small \noindent\textbf{Key words and Phrases:} Random dynamical systems,
Contraction mappings, Perron-Frobenius theory, Nonlinear cocycles, Stochastic
equations, Random monotone mappings, Hilbert-Birkhoff metric.}

{\small \noindent\textbf{2010 Mathematics Subject Classifications:} 37H10,
37H99,\ 37H15.\bigskip}\bigskip\bigskip

\section{Introduction}

Let $V=\mathbb{R}^{n}$ be a finite-dimensional real vector space with some
norm $\Vert\cdot\Vert$. A subset $K$ of $V$ is called a \emph{cone} if it
contains with any vectors $x$ and $y$ any non-negative linear combination
$\alpha x+\beta y$ of these vectors. A cone is called \emph{proper} if
$K\cap(-K)=\{0\}$.

Let $K\subseteq V$ be a closed proper cone in $V$ with non-empty interior
$K^{\circ}$; we will call such cones \emph{solid}. The cone $K$ induces the
partial ordering $\leq_{K}$ in the space $V$ defined as follows: $x\leq_{K}y$
if and only if $y-x\in K$. We shall write $x\prec_{K}y$ if $x\leq_{K}y$,
$x\neq y$, and $x<_{K}y$ if $y-x\in K^{\circ}$.

Let $L$ be another solid cone in $V$. A mapping $D:K\rightarrow L$ is called
\textit{monotone} if $D(x)\leq_{L}D(y)$ for any vectors $x,y\in K$ satisfying
$x\leq_{K}y$. It is called \textit{completely monotone} if each of the
relations $x\leq_{K}y,\ x\prec_{K}y$ or$\ x<_{K}y$ between two vectors $x,y\in
K$ implies the corresponding relation $D(x)\leq_{L}D(y),$ $D(x)\prec_{L}D(y)$
or $D(x)<_{L}D(y)$ between the vectors $D(x),D(y)\in L$. A mapping $D$ is
termed \textit{strictly monotone} if the relation $x\prec_{K}y$ implies
$D(x)<_{L}D(y)$.

Denote by $V^{\ast}$ the dual to the space $V$. Elements of $V^{\ast}$ are
linear functionals $\phi(x)=\langle\phi,x\rangle$ on $V$. For any cone $K$,
denote by
\[
K^{\ast}=\{\phi\in V^{\ast}\ :\ \phi(x)\geq0\ \text{for all}\ x\in K\}
\]
the cone \textit{dual} to $K$. If $K$ is a solid cone, then so is $K^{\ast}%
$(see \cite{Lemmens2},\ Lemma 1.2.4). Every functional in the interior of
$K^{\ast}$ is \textit{strictly positive}, i.e., $\phi(x)>0\ $for all$\ 0\neq
x\in K$. For any linear functional $\phi$ in the interior of $K^{\ast}$, put
\[
\Sigma_{\phi}^{K}=\{x\in K:\ \phi(x)=1\}.
\]
The set $\Sigma_{\phi}^{K}$ is non-empty, compact and convex (\textit{ibid}).

Let $(\Omega,\mathcal{F},P)$ be a complete probability space and
$T:\Omega\rightarrow\Omega$ its automorphism, i.e., a one-to-one mapping such
that $T$ and $T^{-1}$ are measurable and preserve the measure $P$. Let
$...\subseteq\mathcal{F}_{-1}\subseteq\mathcal{F}_{0}\subseteq\mathcal{F}%
_{1}\subseteq...~$be a filtration on $\Omega$ such that each $\sigma$-algebra
$\mathcal{F}_{t}$ is completed by $\mathcal{F}$-measurable sets of measure
$0$. Assume this filtration is invariant with respect to $T$, i.e.
$\mathcal{F}_{t+1}=T^{-1}\mathcal{F}_{t}$ for each $t$. Suppose that for every
$\omega\in\Omega$, a solid cone $K(\omega)\subseteq V$ depending
$\mathcal{F}_{0}$-measurably\footnote{We say that the set-valued mapping
$\omega\mapsto K(\omega)\subseteq V$ is $\mathcal{F}_{0}$-measurable (or the
set $K(\omega)$ depends $\mathcal{F}_{0}$-measurably on $\omega$) if its graph
$\{(\omega,x):$ $x\in K(\omega)\}$ belongs to $\mathcal{F}_{0}\times
\mathcal{V}$, where $\mathcal{V}$ is the Borel $\sigma$-algebra on $V$.} on
$\omega$ is given. Put $K_{t}(\omega)=K(T^{t}\omega)$, $t=0,\pm1,\pm2,...$.
Let $D(\omega,x)$ be a mapping of the cone $K_{0}(\omega)$ into the cone
$K_{1}(\omega)$. Define%
\[
D_{t}(\omega,x)=D(T^{t-1}\omega,x),\ t=0,\pm1,\pm2,....
\]
For shortness, we will write $D(\omega)$ in place of $D(\omega,x)$ and
$D_{t}(\omega)$ in place of $D_{t}(\omega,x)$. Put
\[
C(t,\omega)=D_{t}(\omega)D_{t-1}(\omega)...D_{1}(\omega),\ t=1,2,...,
\]
where the product means the composition of maps, and $C(0,\omega)=Id$ (the
identity map). We have
\[
C(t,T^{s}\omega)C(s,\omega)=C(t+s,\omega),\ t,s\geq0,
\]
i.e., the mapping $C(t,\omega)$ is a \textit{cocycle} over the dynamical
system $(\Omega,\mathcal{F},P,T)$ (see Arnold \cite{Arnold}). The mapping
$C(t,\omega)$ transforms elements of the cone $K_{0}(\omega)$ into elements of
the cone $K_{t}(\omega)$.

It can be shown that the interior of the dual cone $K^{\ast}(\omega)$ depends
$\mathcal{F}_{0}$-measurably on $\omega$, and so there exists an
$\mathcal{F}_{0}$-measurable linear functional $\phi(\omega)$ such that
$\phi(\omega)$ belongs to the interior of $K^{\ast}(\omega)$ for each $\omega$
(see Propositions A.1 and A.3 in the Appendix). We fix the functional
$\phi(\omega)$ and define $\hat{K}(\omega)=\Sigma_{\phi(\omega)}^{K(\omega)}$.
The set $\hat{K}(\omega)$ is non-empty, compact and convex.

Let us extend the mapping $D(\omega,x)$ to all $x\in V$ by setting $\bar
{D}(\omega,x)=D(\omega,x)$ if $x\in K(\omega)$ and $\bar{D}(\omega,x)=\infty$
for $x\notin K(\omega)$, where "$\infty$" stands for a one-point
compactification of $V$. We will impose the following conditions:

(\textbf{D1}) $D(\omega,x)$ is continuous in $x\in K(\omega)$ and $\bar
{D}(\omega,x)$ is $\mathcal{F}_{1}\mathcal{\times V}$-measurable in
$(\omega,x)\in\Omega\times V$, where $\mathcal{V}$ is the Borel $\sigma
$-algebra on $V$.

(\textbf{D2}) $D(\omega,x)$ is positively homogeneous (of degree one) in $x\in
K(\omega)$:
\[
D(\omega,\lambda x)=\lambda D(\omega,x)\ \ \text{for any}\ \lambda>0,\ x\in
K(\omega).
\]

(\textbf{D3}) $D(\omega,x)$ is a completely monotone mapping from
$K_{0}(\omega)$ into $K_{1}(\omega)$.

Furthermore, we will assume that the cocycle $C(t,\omega)$ satisfies the
following condition.

(\textbf{C}) For almost all $\omega\in\Omega$ there is a natural number
$l_{\omega}$ such that the mapping $C(l_{\omega},\omega)$ is strictly monotone.

The main result of this paper is as follows.

\textbf{Theorem 1.}\textit{\ (a) There exist an }$\mathcal{F}_{0}%
$\textit{-measurable vector function }$x(\omega)$\textit{ and an }%
$\mathcal{F}_{1}$\textit{-measurable scalar function }$\alpha(\omega)$\textit{
such that}%
\[
\alpha(\omega)>0,\ x(\omega)\in K^{\circ}(\omega),\ \langle\phi(\omega
),x(\omega)\rangle=1
\]
\textit{for all }$\omega$\textit{ and}%
\begin{equation}
\alpha(\omega)x(T\omega)=D(\omega)x(\omega)\text{ (a.s.).} \label{P-F}%
\end{equation}

\textit{(b) The pair of functions }$(\alpha(\omega),x(\omega))$, \textit{where
}$\alpha(\omega)\geq0$,\textit{ }$x(\omega)\in K(\omega)$\textit{ and
}$\langle\phi(\omega),x(\omega)\rangle=1$, \textit{satisfying (\ref{P-F}) is
determined uniquely up to the equivalence with respect to the measure }%
$P$\textit{. }

\textit{(c) If }$t\rightarrow\infty$\textit{, then }
\begin{equation}
\Vert\frac{C(t,T^{-t}\omega)a}{\langle\phi(\omega),C(t,T^{-t}\omega)a\rangle
}-x(\omega)\Vert\rightarrow0\ \text{\textit{\ (a.s.)}}, \label{lim}%
\end{equation}
\textit{where convergence is uniform in }$0\neq a\in K_{-t}(\omega)$\textit{.}

This result may be regarded as a stochastic nonlinear generalization of the
Perron-Frobenius theorem: $x(\cdot)$ and $\alpha(\cdot)$ play the roles of an
\textquotedblleft eigenvector\textquotedblright\ and an \textquotedblleft
eigenvalue\textquotedblright\ of the random mapping $D(\omega)$ with respect
to the dynamical system $T:\Omega\rightarrow\Omega$. The original versions of
this classical theorem were discovered at the beginning of the twentieth
century by Perron \cite{Perron1}, \cite{Perron2}, who investigated eigenvalues
and eigenvectors of matrices with strictly positive entries, and by Frobenius
\cite{Frobenius1},\cite{Frobenius2},\cite{Frobenius3}, who extended Perron's
results to irreducible nonnegative matrices. Extensions of the
Perron-Frobenius results to nonlinear mappings were obtained by H. Nikaido
\cite{Nikaido}, M. Morishima \cite{Morishima}, T. Fujimoto \cite{Fujimoto}, Y.
Oshime \cite{Oshime1}, \cite{Oshime2}, \cite{Oshime3} and others. Those
extensions were motivated by applications in mathematical economics, in
particular, to the so-called nonlinear Leontief model \cite{S-S}. For reviews
of nonlinear versions of the Perron-Frobenius theory, we refer the reader to
the monographs by Nussbaum \cite{Nussbaum1}, \cite{Nussbaum2} and the papers
by {Kohlberg \cite{Kohlberg} }and Gaubert and Gunawardena \cite{Gaubert}.

The first result on stochastic generalizations of the Perron-Frobenius theorem
for linear maps $D(\omega)$ (non-negative random matrices), was obtained in
\cite{Evstigneev1974}. The result was extended and applied to mathematical
models in statistical physics and evolutionary biology by Arnold et al.
\cite{Arnold1994}. The analysis in \cite{Arnold1994} was based on the use of a
stochastic contraction principle for a suitable metric (the Hilbert-Birkhoff
metric). The paper \cite{Evstigneev1974} employed completely different
methods, suggested by those which are used in the analysis of ergodicity and
mixing properties of Markov chains (cf. Dobrushin \cite{Dobrushin}). The first
stochastic nonlinear analogue of the Perron-Frobenius was obtained in the
paper by Evstigneev and Pirogov \cite{Evstigneev2010}. In that paper,
$D(\omega)$ was a mapping of the set $\mathbb{R}_{+}^{n}$ of non-negative
$n$-dimensional vectors into itself. Now we generalize this result to more
general random cones $K(\omega)\subseteq V$.

Problems related to stochastic (linear and nonlinear) Perron-Frobenius
theorems arise in various areas of pure and applied mathematics, in
particular, in statistical physics, ergodic theory, mathematical biology and
mathematical finance, see, e.g., \cite{Kifer1996}, \cite{Kifer2008},
\cite{Arnold1994}, \cite{Dempster-Evstigneev}. Extensions of this theory to
set-valued mappings $D(\omega,x)$ (\textit{von Neumann-Gale dynamical systems}
\cite{Arnlod-Evst-Gund}, \cite{ESH}) have important applications in
mathematical economics and finance \cite{ESH-growth-model}, \cite{DET}.

Several comments about the assumptions imposed are in order. Let $K$ and $L$
be solid cones in $V$. Consider a \textit{concave} mapping $D:K\rightarrow L$,
i.e. a mapping satisfying
\begin{equation}
D(\theta x+(1-\theta)y)\geq_{L}\theta D(x)+(1-\theta)D(y) \label{concave}%
\end{equation}
for all $x,y\in K$ and $\theta\in\lbrack0,1]$. Clearly, if $D$ is homogeneous,
then $D$ is concave if and only if it is \textit{superadditive}:%

\begin{equation}
D(x+y)\geq_{L}D(x)+D(y). \label{superadditive}%
\end{equation}

For a superadditive mapping $D:K\rightarrow L$, the relation$\ x\prec_{K}y$
between two vectors $x,y\in K$ implies the corresponding relation
$D(x)\prec_{L}D(y)$ between the vectors $D(x),D(y)\in L$ if and only if

\textbf{(M1)} $D(h)\succ_{L}0$ for all $h\succ_{K}0$.

The relation$\ x<_{K}y$ implies the corresponding relation $D(x)<_{L}D(y)$ if
and only if

\textbf{(M2)} $D(h)>_{L}0$ for all $h>_{K}0$.

The mapping $D(x)$ is strictly monotone if and only if

\textbf{(M3)} $D(h)>_{L}0$ for all $h\succ_{K}0$.

We can also see from (\ref{superadditive}) that any superadditive mapping is
monotone. By using this, we obtain that if $D$ is concave and homogeneous,
then \textbf{(M2)} is equivalent to

\textbf{(M4)} $D(h_{\ast})>_{L}0$ for some $h_{\ast}\geq_{K}0$.

Clearly, \textbf{(M4)} follows from \textbf{(M2).} Conversely, \textbf{(M4)}
implies \textbf{(M2)} because for any $h>_{K}0$ we have $h\geq_{K}\lambda
h_{\ast}$, where $\lambda>0$ which yields $D(h)\geq_{L}\lambda D(h_{\ast
})>_{L}0$. Thus, for a concave homogeneous mapping, its complete monotonicity
is equivalent to the validity of \textbf{(M1)} and \textbf{(M2)} (or
\textbf{(M1)} and \textbf{(M4)}), and its strict monotonicity is equivalent to
\textbf{(M3)}.

The paper is organized as follows. In Section 2, a stochastic version of the
fixed point principle which plays key role in the proof of Theorem 1 is
established. The proof of Theorem 1 is given in Section 3. The Appendix
contains statements and short proofs of some general facts regarding
measurable selections which are used in this work.\bigskip

\section{Hilbert-Birkhoff metric}

Given a solid cone $K$ and a strictly positive linear functional $\phi$, the
Hilbert-Birkhoff (\cite{Hilbert},\ \cite{Birkhoff1957}) metric on the set
$Y:=\Sigma_{\phi}^{K}\cap K^{\circ}$ is defined as follows. For any $x,y\in Y$
put%
\[
M(x/y)=\inf\{\beta>0:x\leq_{K}\beta y\},\ m(x/y)=\sup\{\alpha>0:\alpha
y\leq_{K}x\}
\]
and%
\begin{equation}
d(x,y)=\log\left[  \frac{M(x/y)}{m(x/y)}\right]  . \label{d}%
\end{equation}
It can be shown (see \cite{Lemmens2}, Propositions 2.1.1 and 2.5.4) that the
function $d(x,y)$ is a complete separable metric on $Y$,\ and the topology
generated by it on $Y$ coincides with the Euclidean topology on $Y$.
Furthermore, there exists a constant $M>0$ such that
\begin{equation}
\Vert x-y\Vert\leq M(e^{d(x,y)}-1), \label{inequality}%
\end{equation}
for all $x,y\in Y$ (see \cite{Lemmens2}, formula (2.21)).

\textbf{Remark.} An important example of $K$ is the cone $\mathbb{R}_{+}^{n}$
consisting of all non-negative vectors in $V=\mathbb{R}^{n}$. Suppose
$\phi(x)=\sum_{i=1}^{n}x_{i}$ for $x=(x_{1},...,x_{n})$. Then we have
\[
\Sigma_{\phi}^{K}=\{x\geq0:\sum_{i=1}^{n}x_{i}=1\},\ Y=\{x>0:\sum_{i=1}%
^{n}x_{i}=1\},
\]%
\[
M(x/y)=\max_{i}(x_{i}/y_{i}),\ m(x/y)=\min_{j}(x_{j}/y_{j}),
\]%
\[
d(x,y)=\log\left[  \max_{i}(x_{i}/y_{i})\cdot\max_{j}(y_{j}/x_{j})\right]  .
\]
Here, the inequalities $x\geq0$ and $x>0$ are understood coordinate-wise.

Hilbert-Birkhoff metric is a particularly useful tool in the study of monotone
homogeneous maps on cones. A mapping $f:X\rightarrow Y$ from a metric space
$(X,d_{X})$ into a metric space $(Y,d_{Y})$ is called \textit{non-expansive}
if
\[
d_{Y}(f(x),f(y))\leq d_{X}(x,y),\ \text{for all}\ x,y\in X.
\]
It is called \textit{strictly non-expansive} if the inequality in the above
formula is strict for all $x\neq y$ in $X$. The usefulness of Hilbert-Birkhoff
metric lies in the fact that linear, and some nonlinear, mappings of cones are
non-expansive with respect to this metric.

Let $K$ and $L$ be solid cones in $V$ and $\phi_{1}\in(K^{\ast})^{\circ}$,
$\phi_{2}\in(L^{\ast})^{\circ}$. Put $Y_{1}=\Sigma_{\phi_{1}}^{K}\cap
K^{\circ}$, $Y_{2}=\Sigma_{\phi_{2}}^{L}\cap L^{\circ}$ and suppose
$d_{i}(x,y)$ is Hilbert-Birkhoff metric on $Y_{i}$, $i=1,2$.

\textbf{Theorem 2}.\textbf{\ }\textit{If }$f:K\rightarrow L$\textit{\ is a
monotone and homogeneous (of degree 1) mapping such that }$f(x)\succ_{L}%
0$\textit{\ for all }$x\succ_{K}0$\textit{, then the mapping }$g:Y_{1}%
\rightarrow Y_{2}$\textit{\ given by }$g(x)=f(x)/\langle\phi_{2},f(x)\rangle
$\textit{\ is non-expansive with respect to the metric }$d_{1}(x,y)$%
\textit{\ on }$Y_{1}$\textit{\ and the metric }$d_{2}(x,y)$\textit{\ on
}$Y_{2}$\textit{. Moreover, if }$f$\textit{\ is strictly monotone and
homogeneous, then }$g$\textit{\ is strictly non-expansive.}

\textit{Proof.} Let $x,y\in Y_{1}$ and write $\alpha=m(x/y),\ \beta=M(x/y)$.
Since $K$ is a closed cone, we have $\alpha y\leq_{K}x\leq_{K}\beta y$ and so
$\alpha f(y)\leq_{L}f(x)\leq_{L}\beta f(y)$ because $f$ is monotone and
homogeneous. Thus,%

\[
\alpha\frac{\langle\phi_{2},f(y)\rangle}{\langle\phi_{2},f(x)\rangle}%
g(y)\leq_{L}g(x)\leq_{L}\beta\frac{\langle\phi_{2},f(y)\rangle}{\langle
\phi_{2},f(x)\rangle}g(y),
\]
which implies%

\[
d_{2}(g(x),g(y))\leq\log(\frac{\beta}{\alpha})=d_{1}(x,y).
\]

Let $f$ be strictly monotone. If $x,y\in Y_{1}$ and $x\neq y$, we have
$x\neq\lambda y$ for all $\lambda>0$ (otherwise if $x=\lambda y$ for some
$\lambda>0$, then $1=\phi_{1}(y)=\phi_{1}(x)=\lambda\phi_{1}(y)$, which yields
$\lambda=1$ and $x=y$). Then $\alpha y\prec_{K}x\prec_{K}\beta y$, and so
$\alpha f(y)<_{L}f(x)<_{L}\beta f(y)$. Hence there exist $\mu>\alpha$ and
$\tau<\beta$ such that $\mu f(y)\leq_{L}f(x)\leq_{L}\tau f(y)$. So that
\[
d_{2}(g(x),g(y))\leq\log(\frac{\tau}{\mu})<\log(\frac{\beta}{\alpha}%
)=d_{1}(x,y).
\]

The proof is complete.

\section{Stochastic fixed-point principle}

In the proof of Theorem 1, we will use a stochastic generalization of the
following well-known result regarding strictly non-expansive mappings (see,
e.g., \cite{Eisenack}, \cite{Kohlberg}). Let $f$ be a strictly non-expansive
mapping from a compact space $X$ into itself. Then $f$ has a unique fixed
point $\overline{x}$, and $f^{k}(x)\rightarrow\overline{x}$ as $k\rightarrow
\infty$ for each $x\in X$. (We denote by $f^{k}(x)$ the $k$th iterate of the
mapping $f$). A stochastic version of the above contraction principle was
obtained in the paper by Evstigneev and Pirogov \cite{Evstigneev2007}. Here we
establish a more general version of this result. Let us formulate it.

As before, let $(\Omega,\mathcal{F},P)$ be a complete probability space,
$T:\Omega\rightarrow\Omega$ its automorphism, and $...\subseteq\mathcal{F}%
_{-1}\subseteq\mathcal{F}_{0}\subseteq\mathcal{F}_{1}\subseteq...$ a
filtration such that each $\mathcal{F}_{t}$ contains all sets in $\mathcal{F}$
of measure $0$. Let $(V,\mathcal{V})$ be a standard\footnote{A measurable
space is called standard if it is isomorphic to a Borel subset of a complete
separable metric space.} measurable space and let $X(\omega)\subseteq V$ be a
non-empty set depending $\mathcal{F}_{0}$-measurably on $\omega\in\Omega$. Let
$f(\omega,x)$ be a mapping assigning to every $\omega\in\Omega$ and every
$x\in X(\omega)$ an element $f(\omega,x)\in X(T\omega)$. Our main goal in this
section is to provide conditions under which the equation
\begin{equation}
\xi(T\omega)=f(\omega,\xi(\omega))\;\text{(a.s.)} \label{Eq}%
\end{equation}
has a solution in the class of measurable mappings $\xi:\Omega\rightarrow V$
such that $\xi(\omega)\in X(\omega)$ for almost all $\omega$. We also will be
interested in the uniqueness of this solution and properties of its stability.
Equations of the type (\ref{Eq}) arise in connection with various questions of
the theory of random dynamical systems (Arnold 1998 \cite{Arnold}). Our study
of such equations is motivated by their applications in the stochastic
Perron-Frobenius theory (Evstigneev 1974 \cite{Evstigneev1974}, Arnold,
Demetrius and Gundlach 1994 \cite{Arnold1994}, Kifer 1996 \cite{Kifer1996} and
Evstigneev and Pirogov 2010 \cite{Evstigneev2010}).

Let us extend $f(\omega,x)$ to the whole space $V$ by setting $\bar{f}%
(\omega,x)=f(\omega,x)$ if $x\in X(\omega)$ and $\bar{f}(\omega,x)=\infty$ if
$x\notin X(\omega)$, where the symbol "$\infty$" denotes a point added to $V$.

Assume that the following conditions hold.

(\textbf{A1}) The mapping $\bar{f}(\omega,x)$ ($\omega\in\Omega$, $x\in V$) is
$\mathcal{F}_{1}\mathcal{\times V}$-measurable.

For each $\omega$, let $Y(\omega)$ be a non-empty subset of $X(\omega)$
equipped with a separable metric $\rho(\omega,x,y)$, $x,y\in Y(\omega)$. Let
us introduce the following assumptions.

(\textbf{A2}) (a) The set-valued mapping $\omega\mapsto Y(\omega)$ is
$\mathcal{F}_{0}$-measurable.

(b) The function
\[
\bar{\rho}(\omega,x,y):=\left\{
\begin{array}
[c]{cc}%
\rho(\omega,x,y), & \text{if }x,y\in Y(\omega),\\
+\infty, & \text{otherwise,}%
\end{array}
\right.
\]
is $\mathcal{F}_{0}\times\mathcal{V}\times\mathcal{V}$-measurable.

(c) For each $\omega$, the Borel measurable structure on $Y(\omega)$ induced
by the metric $\rho(\omega,x,y)$ coincides with the measurable structure
induced on $Y(\omega)$ by the $\sigma$-algebra $\mathcal{V}$.

(d) For each $\omega\in\Omega$ and $x\in Y(\omega)$, we have $f(\omega,x)\in
Y(T\omega)$, and the mapping $f(\omega,\cdot):$ $Y(\omega)\rightarrow
Y(T\omega)$ is continuous with respect to the metric $\rho(\omega,x,y)$ on
$Y(\omega)$ and the metric $\rho(T\omega,x,y)$ on $Y(T\omega)$.

Note that $Y(\omega)\in\mathcal{V}$ for each $\omega$ by virtue of (a) and
that (a) follows from (b).

For every $k=0,\pm1,\pm2,...$ define%
\[
X_{k}(\omega)=X(T^{k}\omega),\ Y_{k}(\omega)=Y(T^{k}\omega),\ \rho_{k}%
(\omega,x,y)=\rho(T^{k}\omega,x,y),
\]%
\begin{equation}
f_{k}(\omega,x)=f(T^{k-1}\omega,x)\ [x\in X_{k-1}(\omega)]. \label{f1}%
\end{equation}
For each $m=0,1,2,...$ put%
\begin{equation}
f^{(m)}(\omega,x)=f_{0}(\omega)f_{-1}(\omega)...f_{-m}(\omega)(x)\ [x\in
X_{-m-1}(\omega)], \label{f2}%
\end{equation}%
\[
X^{(m)}(\omega)=f^{(m)}(\omega,X_{-m-1}(\omega)).
\]
The product $f^{(m)}(\omega,x)=f_{0}(\omega)f_{-1}(\omega)...f_{-m}(\omega)$
means the composition of the mappings. Note that for each $m=0,1,...$ the map
$f_{-m}(\omega,x)$ acts from $\Omega\times X_{-m-1}(\omega)$\ into
$X_{-m}(\omega)$, and so $f^{(m)}(\omega,x)$ acts from $\Omega\times
X_{-m-1}(\omega)$\ into $X_{0}(\omega)$. The extended mappings $\bar{f}%
_{-m}(\omega,x)$\ and $\bar{f}^{(m)}(\omega,x)$ are $\mathcal{F}_{-m}%
\times\mathcal{V}$-measurable and $\mathcal{F}_{0}\times\mathcal{V}%
$-measurable, respectively. The functions $\bar{\rho}_{-m}(\omega,x,y)$
are\textit{\ }measurable with respect to $\mathcal{F}_{-m}\times
\mathcal{V}\times\mathcal{V}\subseteq\mathcal{F}_{0}\times\mathcal{V}%
\times\mathcal{V}$.

(\textbf{A3}) There is a sequence of $\mathcal{F}_{0}$-measurable sets
$\Omega_{0}\subseteq\Omega_{1}\subseteq...\subseteq\Omega$ such that
$P(\Omega_{m})\rightarrow1$ and for each $m=0,1,...$ and $\omega\in\Omega_{m}$
the following conditions are satisfied:

(a) the set $X^{(m)}(\omega)$ is contained in $Y(\omega)$ and is compact with
respect to the metric $\rho(\omega,x,y)$;

(b) for all $x,y\in Y_{-m-1}(\omega)$ with $x\neq y$, we have
\begin{equation}
\rho(\omega,f^{(m)}(\omega,x),f^{(m)}(\omega,y))<\rho_{-m-1}(\omega,x,y).
\label{non-expansive}%
\end{equation}
Since the sequence of sets $\Omega_{m}$ is non-decreasing, there exists an
$\mathcal{F}_{0}$-measurable function $m(\omega)$ with non-negative integer
values such that for each $\omega\in\bar{\Omega}:=\Omega_{1}\cup\Omega_{2}%
\cup...$ (and hence for almost all $\omega$), we have $\omega\in\Omega_{m}%
$,$\ m\geq m(\omega)$. We can define $m(\omega)=\min\{i:\omega\in\Omega_{i}\}$
if $\omega\in\bar{\Omega}$ and $m(\omega)=0$, otherwise.

\textbf{Theorem 3.}\textit{\ (i) There exists an }$\mathcal{F}_{0}%
$-\textit{measurable mapping }$\xi:\Omega\rightarrow V$ \textit{such that
}$\xi(\omega)\in Y(\omega)$\textit{, equation (\ref{Eq}) holds, and}
\begin{equation}
\lim_{m(\omega)\leq m\rightarrow\infty}\sup_{x\in X_{-m-1}(\omega)}\rho
(\omega,\xi(\omega),f_{0}(\omega)...f_{-m}(\omega
)(x))=0\ \text{\textit{(a.s.)}}. \label{x-converg}%
\end{equation}
\textit{\ }

\textit{(ii) If }$\eta:\Omega\rightarrow V$\textit{\ is any (not necessarily
measurable) mapping for which }$\eta(\omega)\in X(\omega)$ \textit{and
equation (\ref{Eq}) holds,\ then }$\eta=\xi$\textit{\ with probability one.}

According to (\ref{x-converg}), the sequence $f_{0}...f_{-m}(x)$ converges to
$\xi(\omega)$ in the metric $\rho(\omega,x,y)$ uniformly in $x\in
X_{-m-1}(\omega)$ with probability one. Note that the distance $\rho
(\omega,\cdot,\cdot)$ between $f_{0}...f_{-m}(x)$ and $\xi(\omega)$ involved
in (\ref{x-converg}) is defined only if $f_{0}...f_{-m}(x)\in Y(\omega)$. By
virtue of condition (a) in (\textbf{A3}), this inclusion holds for almost all
$\omega$, all $m\geq m(\omega)$ and $x\in X_{-m-1}(\omega)$, therefore the
limit in (\ref{x-converg}) is taken over $m\geq m(\omega)$.

\section{Proof of the stochastic fixed point principle}

\textit{Proof of Theorem 3.} \textit{1st step.} Observe that $X^{(0)}%
(\omega)\supseteq X^{(1)}(\omega)\supseteq X^{(2)}(\omega)\supseteq...$ and
$X^{(m)}(\omega)\neq\emptyset$ for each $m$ and $\omega$. Consider the sets
$\Omega_{m}$ ($m=0,1,...$) described in (\textbf{A3}) and their union
$\bar{\Omega}$ . According to (\textbf{A3}), $P(\bar{\Omega})=1$ and each
$\omega\in\bar{\Omega}$ belongs to all $\Omega_{m},\,m\geq m(\omega)$. For
$\omega\in$ $\bar{\Omega}$, all the sets $X^{(m)}(\omega),\;m\geq m(\omega)$,
are contained in $Y(\omega)$ and compact, and so the set $X^{\infty}%
(\omega):=\cap_{m=0}^{\infty}X^{(m)}(\omega)\subseteq Y(\omega)$ is non-empty
and compact as an intersection of a nested sequence of non-empty compacta
$X^{(m)}(\omega),\;m\geq m(\omega)$.

\textit{2nd step. }Define $\Omega^{\ast}=\cap_{k=-\infty}^{+\infty}(T^{k}%
\bar{\Omega})$. The set $\Omega^{\ast}$ is invariant and $P(\Omega^{\ast})=1$.
Let us show that
\begin{equation}
X^{\infty}(T\omega)=f(\omega,X^{\infty}(\omega)),\;\omega\in\Omega^{\ast}.
\label{e1}%
\end{equation}
Equality (\ref{e1}) is equivalent to
\begin{equation}
X^{\infty}(\omega)=f(T^{-1}\omega,X^{\infty}(T^{-1}\omega)),\;\omega\in
\Omega^{\ast}, \label{e2}%
\end{equation}
because $\omega\in\Omega^{\ast}$ if and only if $T^{-1}\omega\in\Omega^{\ast}%
$. To prove (\ref{e2}) let us observe that
\begin{equation}
f(T^{-1}\omega,%
{\textstyle\bigcap\limits_{m=0}^{\infty}}
X^{(m)}(T^{-1}\omega))=%
{\textstyle\bigcap\limits_{m=0}^{\infty}}
f(T^{-1}\omega,X^{(m)}(T^{-1}\omega)),\;\omega\in\Omega^{\ast}. \label{e4}%
\end{equation}
The inclusion \textquotedblright$\subseteq$\textquotedblright\ in (\ref{e4})
holds always. The opposite inclusion follows from the continuity of
$f(T^{-1}\omega,\cdot)$ on $Y(T^{-1}\omega)$ and the fact that $X^{(m)}%
(T^{-1}\omega)$ are nested and compact in $Y(T^{-1}\omega)$ for all $m$ large
enough. By using (\ref{e4}), we obtain
\[
f(T^{-1}\omega,X^{\infty}(T^{-1}\omega))=f(T^{-1}\omega,%
{\textstyle\bigcap\limits_{m=0}^{\infty}}
X^{(m)}(T^{-1}\omega))
\]%
\[
=%
{\textstyle\bigcap\limits_{m=0}^{\infty}}
f(T^{-1}\omega,X^{(m)}(T^{-1}\omega))=%
{\textstyle\bigcap\limits_{m=0}^{\infty}}
f_{0}(\omega,X^{(m)}(T^{-1}\omega))
\]%
\[
=%
{\textstyle\bigcap\limits_{m=0}^{\infty}}
X^{(m+1)}(\omega)=X^{\infty}(\omega),\;\omega\in\Omega^{\ast}.
\]
The fourth equality in this chain of relations holds because
\[
X^{(m)}(T^{-1}\omega)=f_{0}(T^{-1}\omega)f_{-1}(T^{-1}\omega)...f_{-m}%
(T^{-1}\omega)(X_{-m-1}(T^{-1}\omega))
\]%
\[
=f_{-1}(\omega)f_{-2}(\omega)...f_{-m-1}(\omega)(X_{-m-2}(\omega)),
\]
and so
\[
f_{0}(\omega)(X^{(m)}(T^{-1}\omega))=f_{0}(\omega)f_{-1}(\omega)...f_{-m-1}%
(\omega)(X_{-m-2}(\omega))=X^{(m+1)}(\omega).
\]

\textit{3rd step.} For $\omega\in\Omega^{\ast}$, denote the diameter in the
metric $\rho(\omega,x,y)$ of the compact set $X^{\infty}(\omega)\subseteq
Y(\omega)$ by $\rho(\omega)$ and put $\rho(\omega)=+\infty$ if $\omega
\in\Omega\setminus\Omega^{\ast}$. For $m=0,1,2,...$, put $\Omega_{m}^{\ast
}:=\Omega^{\ast}\cap\Omega_{m}$ and for $\omega\in\Omega$ define
\begin{equation}
\rho^{(m)}(\omega)=\left\{
\begin{array}
[c]{cc}%
\text{\textrm{diam}\thinspace}X^{(m)}(\omega), & \text{if }\omega\in\Omega
_{m}^{\ast},\\
+\infty, & \text{otherwise.}%
\end{array}
\right.  \label{ro-k}%
\end{equation}
Recall that, for $\omega\in\Omega_{m}$ and hence for $\omega\in\Omega
_{m}^{\ast}$, the set $X^{(m)}(\omega)$ is contained in $Y(\omega)$ and is
compact, so that its diameter \textrm{diam}\thinspace$X^{(m)}(\omega)$ in the
metric $\rho(\omega,x,y)$ is well-defined and finite. We claim that
$\rho^{(m)}(\omega)$ is an $\mathcal{F}_{0}$-measurable function of $\omega
\in\Omega$. To prove this assertion we observe that for $\omega\in\Omega
_{m}^{\ast}$, we have \textrm{diam}\thinspace$X^{(m)}(\omega)=$
\textrm{diam\thinspace}$f^{(m)}(\omega,X_{-m-1}(\omega))$,$\;$where
\[
f^{(m)}(\omega,x)=f_{0}(\omega)f_{-1}(\omega)...f_{-m}(\omega)(x),x\in
X_{-m-1}(\omega).
\]
Consequently, for each real $a$, the set $\Omega_{m}^{a}$ of $\omega\in
\Omega_{m}^{\ast}$ satisfying \textrm{diam}\thinspace$X^{(m)}(\omega)>a$ is
the projection on $\Omega_{m}^{\ast}$ of the set
\begin{equation}
\{(\omega,x,y)\in\Omega_{m}^{\ast}\times X_{-m-1}(\omega)\times X_{-m-1}%
(\omega):\;\rho(\omega,f^{(m)}(\omega,x),\,f^{(m)}(\omega,y))>a\text{ }\},
\label{ro-meas}%
\end{equation}
which is an $\mathcal{F}_{0}\times\mathcal{V}\times\mathcal{V}$-measurable
subset in $\Omega_{m}^{\ast}\times V\times V$ by virtue of assumptions
(\textbf{A1}) and (\textbf{A2}). Since $V$ (and hence $V\times V$) is standard
and $(\Omega,\mathcal{F}_{0},P)$ is a complete probability space, $\Omega
_{m}^{a}$ is $\mathcal{F}_{0}$-measurable (see, e.g., Dellacherie and Meyer
1978, Theorem III.33). This implies that $\rho^{(m)}(\omega)$ is
$\mathcal{F}_{0}$-measurable because $\rho^{(m)}(\omega)=+\infty$ outside
$\Omega_{m}^{\ast}$. Finally, $\rho(\omega)$ is $\mathcal{F}_{0}$-measurable
because
\begin{equation}
\rho(\omega)=\lim_{m\rightarrow\infty}\rho^{(m)}(\omega)\text{ for }\omega
\in\Omega^{\ast}, \label{lim-ro-k}%
\end{equation}
which follows the fact that $X^{(m)}(\omega)$ are nested and compact in
$Y(\omega)$ for all $\omega\in\Omega^{\ast}$ and $m\geq m(\omega)$.

\textit{4th step.} Let us show that $\rho(\omega)=0$ (a.s.). Observe that
equality (\ref{e1}) implies%
\[
X^{\infty}(\omega)=f(T^{-1}\omega,X^{\infty}(T^{-1}\omega))=f(T^{-1}%
\omega)(X^{\infty}(T^{-1}\omega))
\]%
\[
=f(T^{-1}\omega)f(T^{-2}\omega)(X^{\infty}(T^{-2}\omega))=...=f(T^{-1}%
\omega)...f(T^{-m-1}\omega)(X^{\infty}(T^{-m-1}\omega))
\]%
\begin{equation}
=f_{0}(\omega)...f_{-m}(\omega)(X^{\infty}(T^{-m-1}\omega))=f^{(m)}%
(\omega,X^{\infty}(T^{-m-1}\omega)),\;\omega\in\Omega^{\ast}. \label{X-m}%
\end{equation}
By virtue of (\ref{X-m}) and condition (b) in \textbf{(A3)}, for $\omega
\in\Omega_{m}^{\ast}$, we have
\begin{equation}
\rho(\omega)=\rho(\omega,X^{\infty}(\omega))\leq\rho(T^{-m-1}\omega,X^{\infty
}(T^{-m-1}\omega))=\rho(T^{-m-1}\omega) \label{ro}%
\end{equation}
and\qquad\
\begin{equation}
\text{if }\rho(\omega)>0\text{, then }\rho(\omega)<\rho(T^{-m-1}\omega).
\label{ro1}%
\end{equation}
(We also use here the fact that $X^{\infty}(\omega)$ is compact.) Since
$P(\Omega_{m}^{\ast})=P(\Omega^{\ast}\cap\Omega_{m})\rightarrow1$, inequality
(\ref{ro}) yields
\begin{equation}
\lim_{m\rightarrow\infty}P\{\rho(\omega)\leq\rho(T^{-m}\omega)\}\rightarrow1.
\label{ro3}%
\end{equation}

We claim that (\ref{ro3}) implies
\begin{equation}
\rho(\omega)=\rho(T^{-m}\omega)\;\text{a.s. for all }m. \label{ro-inv}%
\end{equation}
To deduce (\ref{ro-inv}) from (\ref{ro3}) we may assume that $\rho(\omega)$ is
bounded by some constant $C$ (we can always replace $\rho(\omega)$ by
$\arctan\rho(\omega)$). By setting $\Delta_{m}:=\{\omega:\rho(\omega)\leq
\rho(T^{-m}\omega)\}$, we write
\[
E|\rho(\omega)-\rho(T^{-m}\omega)|\leq E(\rho(T^{-m}\omega)-\rho(\omega
))\chi_{\Delta_{m}}+CP(\Omega\setminus\Delta_{m}),
\]
where $\chi_{\Delta_{m}}$ is the indicator function of $\Delta_{m}$. Further,
since $E\rho(T^{-m}\omega)=E\rho(\omega)$, we have
\begin{align*}
E(\rho(T^{-m}\omega)-\rho(\omega))\chi_{\Delta_{m}}  &  =E(\rho(T^{-m}%
\omega)-\rho(\omega))\chi_{\Delta_{m}}-E(\rho(T^{-m}\omega)-\rho(\omega))\\
&  =-E(\rho(T^{-m}\omega)-\rho(\omega))\chi_{\Omega\setminus\Delta_{m}}\leq
CP(\Omega\setminus\Delta_{m}).
\end{align*}
Consequently,
\[
E|\rho(\omega)-\rho(T^{-m}\omega)|\leq2CP(\Omega\setminus\Delta_{m}%
)\rightarrow0,
\]
which implies (\ref{ro-inv}).

Suppose $\rho(\omega)>0$ with strictly positive probability. Then there exists
a number $m$ and a set $\Gamma\in\mathcal{F}_{0}$ contained in $\Omega
_{m}^{\ast}$ such that $P(\Gamma)>0$ and $\rho(\omega)>0$ on $\Gamma$. By
virtue of (\ref{ro1}), we have $\rho(\omega)<\rho(T^{-m-1}\omega)\ $for
$\omega\in\Gamma$. On the other hand, we proved that $\rho(\omega
)=\rho(T^{-m-1}\omega)$ for almost all $\omega$. A contradiction.

\textit{5th step.} Since the $\mathcal{F}_{0}$-measurable function
$\rho(\omega)$ is zero a.s., there is a set $\tilde{\Omega}\in\mathcal{F}_{0}$
of full measure such that
\begin{equation}
\tilde{\Omega}\subseteq\Omega^{\ast}\text{ and }\rho(\omega)=0\;\text{for each
}\omega\in\tilde{\Omega}. \label{omega-tilde}%
\end{equation}
This means that for $\omega\in\tilde{\Omega}$, the set $X^{\infty}(\omega
)\ $consists of exactly one point, $\xi^{\infty}(\omega)$. Replacing
$\tilde{\Omega}$ by $\cap_{k=-\infty}^{+\infty}(T^{k}\tilde{\Omega})$, we may
assume that $\tilde{\Omega}$ is invariant.

For every $\omega$, fix any point $\tilde{y}(\omega)$ in the non-empty set
$Y(\omega)$ and put $\xi(\omega)=\xi^{\infty}(\omega)$ for $\omega\in
\tilde{\Omega}$ and $\xi(\omega)=\tilde{y}(\omega)$ for $\omega\in
\Omega\setminus\tilde{\Omega}$. Then for any $\omega\in\tilde{\Omega}%
\subseteq\Omega^{\ast}$ we have $T\omega\in\tilde{\Omega}\subseteq\Omega
^{\ast}$, and so
\begin{align*}
\{\xi(T\omega)\}  &  =\{\xi^{\infty}(T\omega)\}=X^{\infty}(T\omega
)=f(\omega,X^{\infty}(\omega))\\
&  =f(\omega,\{\xi^{\infty}(\omega)\})=f(\omega,\{\xi(\omega)\}\}
\end{align*}
by virtue of (\ref{e1}). Consequently, $\xi(\omega)$ satisfies (\ref{Eq}) for
all $\omega$ in the set $\tilde{\Omega}\subseteq\Omega^{\ast}\subseteq
\bar{\Omega}$ of measure one.

Consider the functions $\rho^{(m)}(\omega)$ defined by (\ref{ro-k}). For each
$\omega\in\tilde{\Omega}$ and $m\geq m(\omega)$ we have $\omega\in\Omega_{m}$,
and so
\begin{equation}
\sup_{x\in X_{-m-1}(\omega)}\rho(\omega,\xi(\omega),f_{0}...f_{-m}%
(x))\leq\text{\textrm{diam\thinspace}}X^{(m)}(\omega)=\rho^{(m)}(\omega)
\label{ro-diam}%
\end{equation}
because $\{\xi(\omega)\}=X^{\infty}(\omega)\subseteq X^{(m)}(\omega)$. This
implies (\ref{x-converg}) since $\lim\rho^{(m)}(\omega)=\rho(\omega)=0$ on the
set $\tilde{\Omega}$ of full measure.

\textit{6th step.} To complete the proof of (i) it is sufficient to show that
the mapping $\xi$ constructed above coincides a.s. with some $\mathcal{F}_{0}%
$-measurable mapping $\zeta$. Then $\zeta$ will be the sought-for solution to
(\ref{Eq}) possessing the properties listed in (i).

Consider some $\mathcal{F}_{0}$-measurable mappings $y_{0}(\omega)$,
$x_{-m-1}(\omega)$, $m=0,1,...$,\ with values in $V$ such that $y_{0}%
(\omega)\in Y(\omega)$, $x_{-m-1}(\omega)\in X_{-m-1}(\omega)$ for all
$\omega$. The existence of these mappings follows from the measurable
selection theorem (see Proposition A.1 in the Appendix) because the graphs of
the set-valued mappings $\omega\mapsto Y(\omega)$ and $X_{-m-1}(\omega)$ are
measurable with respect to $\mathcal{F}_{0}\times\mathcal{V}$ and
$\mathcal{F}_{-m-1}\times\mathcal{V\subseteq F}_{0}\times\mathcal{V}$,
respectively. Define the mappings $\zeta^{m}(\omega)$ of $\Omega$ into $V$
($m=0,1,...$) by the formula
\[
\zeta^{m}(\omega)=\left\{
\begin{array}
[c]{cc}%
f_{0}(\omega)...f_{-m}(\omega)(x_{-m-1}(\omega)) & \text{if }\omega\in
\Omega_{m}\text{,}\\
y_{0}(\omega)\,, & \text{otherwise.}%
\end{array}
\right.
\]
Clearly $\zeta^{m}(\omega)\in Y(\omega)$ for all $\omega$. The mappings
$\zeta^{m}(\omega)$ are $\mathcal{F}_{0}$-measurable because $\Omega_{m}%
\in\mathcal{F}_{0}$ and the mappings $\bar{f}_{-m}(\omega,x)$ are measurable
with respect to $\mathcal{F}_{-m-1}\times\mathcal{V\subseteq}$ $\mathcal{F}%
_{0}\times\mathcal{V}$. For each $\omega\in\tilde{\Omega}$ and $m\geq
m(\omega)$, we have $\omega\in\Omega_{m}$ and
\[
\rho(\omega,\xi(\omega),\zeta^{m}(\omega))=\rho(\omega,\xi(\omega
),f_{0}(\omega)...f_{-m}(\omega)(x_{-m-1}(\omega))\leq\rho^{(m)}(\omega)
\]
where $\rho^{(m)}(\omega)\rightarrow0$ as $m\rightarrow\infty$ (see
(\ref{lim-ro-k}), (\ref{omega-tilde}) and (\ref{ro-diam})). Thus $\zeta
^{m}(\omega)\rightarrow\xi(\omega)$ on a set $\tilde{\Omega}\in\mathcal{F}%
_{0}$, where $P(\tilde{\Omega})=1$. Thus $\xi(\omega)$ is an a.s. limit of
$\mathcal{F}_{0}$-measurable functions, and consequently, it is $\mathcal{F}%
_{0}$-measurable since $\mathcal{F}_{0}$ is complete.

We know that $\xi(T\omega)=f(\omega,\xi(\omega))$ (a.s.) and since $\xi
(\omega)$ coincides a.s. with some $\mathcal{F}_{0}$-measurable mapping
$\xi^{\prime}(\omega)$, we obtain%
\[
\xi^{\prime}(T\omega)=\xi(T\omega)=f(\omega,\xi(\omega))=f(\omega,\xi^{\prime
}(\omega))\text{ (a.s.)},
\]
where the first equality is valid because the transformation $T$ preserves the
measure $P$.

\textit{7th step.} It remains to prove (ii). If $\eta:\Omega\rightarrow V$ is
a mapping for which $\eta(\omega)\in X(\omega)$ and equation (\ref{Eq}) holds
a.s., then%
\[
\eta(\omega)=f(T^{-1}\omega,\eta(T^{-1}\omega))=f(T^{-1}\omega)(\eta
(T^{-1}\omega))
\]
\[
=f(T^{-1}\omega)f(T^{-2}\omega)(\eta(T^{-2}\omega))=...=f(T^{-1}%
\omega)...f(T^{-m-1}\omega)(\eta(T^{-m-1}\omega))\;\text{(a.s.),}%
\]
which yields
\begin{equation}
\eta(\omega)=f_{0}(\omega)f_{-1}(\omega)...f_{-m}(\omega)(\eta(T^{-m-1}%
\omega))\;\text{(a.s.).} \label{et}%
\end{equation}
By combining (\ref{et}) and (\ref{x-converg}), we get
\[
\rho(\omega,\xi(\omega),\eta(\omega))\leq\sup_{x\in X_{-m-1}(\omega)}%
\rho(\omega,\xi(\omega),f_{0}(\omega)...f_{-m}(\omega)(x))\rightarrow
0\ \;\text{(a.s.),}%
\]
and so $\xi(\omega)=\eta(\omega)$ (a.s.). The proof is complete.\bigskip

\section{Nonlinear Perron-Frobenius Theorem: proof}

In this section we prove Theorem 1. The proof is based on a lemma.

\textbf{Lemma 1.}\textit{\ There exists a sequence of }$\mathcal{F}_{0}%
$\textit{-measurable sets }$\Gamma_{1}\subseteq\Gamma_{2}\subseteq
...\subseteq\Omega$\textit{\ such that }$P(\Gamma_{m})\rightarrow
1$\textit{\ and for each }$m=1,2,...\text{ and}\ \omega\in\Gamma_{m}$\textit{,
the mapping }$C(m,T^{-m}\omega)$\textit{\ from the cone }$K_{-m}(\omega
)$\textit{\ to the cone }$K_{0}(\omega)$ \textit{is strictly monotone.}

\textit{Proof. }For each $m\geq1$, consider the set $\Delta_{m}$ of those
$\omega$ for which the mapping $C(m,\omega,x)=C(m,\omega)x$ of the cone
$K_{0}(\omega)$ into the cone $K_{m}(\omega)$ is strictly monotone in $x$. Let
us show that $\Delta_{m}\in\mathcal{F}_{m}$. Denote by $H_{m}(\omega)$ the
closed set $V\setminus K_{m}^{\circ}(\omega)$ and by $\delta(z,H_{m}(\omega))$
the distance (defined in terms of the norm $\Vert\cdot\Vert$) between the
point $z\in V$ and $H_{m}(\omega)$. Clearly, $z\in K_{m}^{\circ}(\omega)$ if
and only if $\delta(z,H_{m}(\omega))>0$.

For each $i,j=1,2,...$ denote by $\Lambda_{ij}(\omega)$ the set of those
$(x,y)\in K_{0}(\omega)\times K_{0}(\omega)$ for which%
\[
y-x\in K_{0}(\omega),\ \Vert y-x\Vert\geq1/i,\ \max\{\Vert x\Vert,\Vert
y\Vert\}\leq j.
\]
The set-valued mapping $\omega\mapsto\Lambda_{ij}(\omega)$ is $\mathcal{F}%
_{0}$-measurable, and so it possesses a countable dense set of $\mathcal{F}%
_{0}$-measurable selections $(x_{ij}^{l}(\omega),y_{ij}^{l}(\omega))$,
$l=1,2,...$(see Proposition A.2 in the Appendix). In view of the compactness
of $\Lambda_{ij}(\omega)$ and continuity of $C(m,\omega,\cdot)$, we have
\[
\Delta_{m}=\bigcup\limits_{i,j=1}^{\infty}\{\omega:\inf\limits_{l}%
[\delta(C(m,\omega,y_{ij}^{l}(\omega))-C(m,\omega,x_{ij}^{l}(\omega
)),H_{m}(\omega))]>0\}.
\]
Since the set-valued mapping $\omega\mapsto H_{m}(\omega)$ is $\mathcal{F}%
_{m}$-measurable, the set $\Delta_{m}$ is a union of a countable family of
$\mathcal{F}_{m}$-measurable sets and is thus $\mathcal{F}_{m}$-measurable
(see Propositions A.2 and A.3 in the Appendix).

If $\omega\in\Delta_{m}$ and $y\succ_{K_{0}(\omega)}x$, we have $C(m,\omega
)y>_{K_{m}(\omega)}C(m,\omega)x$. Furthermore, $D(T^{m}\omega)$ is a
completely monotone mapping from $K_{0}(T^{m}\omega)=K_{m}(\omega)$ into
$K_{1}(T^{m}\omega)=K_{m+1}(\omega)$. Therefore
\[
C(m+1,\omega)y=D(T^{m}\omega)C(m,\omega)y>_{K_{m+1}(\omega)}D(T^{m}%
\omega)C(m,\omega)x=C(m+1,\omega)x
\]
and so $\omega\in\Delta_{m+1}$. Consequently, $\Delta_{m}\subseteq\Delta
_{m+1}$. By virtue of assumption (\textbf{C}), we have $P(\bigcup
_{m=1}^{\infty}\Delta_{m})=1$. By virtue of the inclusion $\Delta_{m}%
\subseteq\Delta_{m+1}$, this implies $P(\Delta_{m})\rightarrow1$. Define
$\Gamma_{m}=T^{m}\Delta_{m}$. Then $\omega\in\Gamma_{m}$ if and only if the
mapping $C(m,T^{-m}\omega)$ of the cone $K_{0}(T^{-m}\omega)=K_{-m}(\omega)$
into the cone $K_{m}(T^{-m}\omega)=K_{0}(\omega)$ is strictly monotone.
Furthermore, $\Gamma_{m}\in\mathcal{F}_{0}$ for every $m$ because $T^{m}%
\Delta_{m}\in\mathcal{F}_{0}$ if and only if $\Delta_{m}\in\mathcal{F}_{m}$.
If $\omega\in\Gamma_{m}$, then $\omega\in\Gamma_{m+1}$ because the mapping
\begin{equation}
C(m+1,T^{-m-1}\omega)=C(m,T^{-m}\omega)D(T^{-m-1}\omega) \label{C-D}%
\end{equation}
of $K_{-m-1}(\omega)$ into $K_{0}(\omega)$ is strictly monotone as the product
of two mappings one of which is completely monotone (from $K_{-m-1}(\omega)$
to $K_{-m}(\omega)$) and the other strictly monotone\ (from $K_{-m}(\omega)$
to $K_{0}(\omega)$). Thus, $\Gamma_{1}\subseteq\Gamma_{2}\subseteq\ldots$,
where $P(\Gamma_{m})=P(\Delta_{m})\rightarrow1$, which completes the proof.

\textit{Proof of Theorem 1.} Put
\[
X(\omega)=\hat{K}(\omega),\ \phi_{t}(\omega)=\phi(T^{t}\omega),\ \hat{K}%
_{t}(\omega)=\hat{K}(T^{t}\omega),\ X_{t}(\omega)=X(T^{t}\omega)\
\]
for any $t=0,\pm1,\pm2,...$. We will apply Theorem 3 to the mapping
\begin{equation}
f(\omega,x)=\frac{D(\omega,x)}{\langle\phi_{1}(\omega),D(\omega,x)\rangle
},\ x\in X_{0}(\omega). \label{f}%
\end{equation}
The mapping $f(\omega,x)$ is well-defined because $D(\omega,x)$ is a
completely monotone mapping from $K_{0}(\omega)\ $into $K_{1}(\omega)$. This
implies that $\langle\phi_{1}(\omega),D(\omega,x)\rangle>0$ because
$D(\omega,x)\neq0$ for $x\neq0$. Furthermore, $\langle\phi_{1}(\omega
),f(\omega,x)\rangle=1$, which means that $f(\omega,x)\in X_{1}(\omega
)=\hat{K}_{1}(\omega)$, and so $f(\omega,\cdot)$ is a mapping of $X_{0}%
(\omega)$ into $X_{1}(\omega)$.

Denote by $\mathcal{V}$ the Borel $\sigma$-algebra on $V$ induced by the
Euclidean topology on $V$. Then the measurable space $(V,\mathcal{V})$ is
standard. Let us verify the assumptions of Theorem 3 for the mappings
$\omega\mapsto X(\omega)$ and $f(\omega,x):X_{0}(\omega)\rightarrow
X_{1}(\omega)$.

We have
\[
\{(\omega,x):x\in X(\omega)\}=\{(\omega,x):x\in K(\omega),\ \langle\phi
(\omega),x\rangle=1\}\in\mathcal{F}_{0}\times\mathcal{V}%
\]
because $\phi(\omega)$ and $\omega\mapsto K(\omega)$ are $\mathcal{F}_{0} $-measurable.

To check (\textbf{A1}) we need to show that the mapping $f(\omega,x)$, which
is equal to $f(\omega,x)$ if $x\in X_{0}(\omega)$ and $\infty$ otherwise, is
$\mathcal{F}_{1}\times\mathcal{V}$-measurable. This follows from the fact that
the set $\Gamma:=\{(\omega,x):$ $\omega\in\Omega,\ $\ $x\in K_{0}(\omega
),\ $\ $\langle\phi_{0}(\omega),x\rangle=1\}$ is $\mathcal{F}_{0}%
\times\mathcal{V}$-measurable and the mapping $f(\omega,x)$ (see (\ref{f}))
restricted to $\Gamma$ is $\mathcal{F}_{1}\times\mathcal{V}$-measurable by
virtue of (\textbf{D1}) and $\mathcal{F}_{1}$-measurability of $\phi
_{1}(\omega)$.

For each $\omega$, we define $Y(\omega)$ as $\hat{K}(\omega)\cap K^{\circ
}(\omega)$ (which corresponds to our previous notation) and consider the
Hilbert-Birkhoff metric $d(\omega)=d(\omega,x,y)$ on $Y(\omega)$. For every
$k=0,\pm1,\pm2,...$ define%

\[
Y_{t}(\omega)=Y(T^{t}\omega),\ d_{k}(\omega)=d_{k}(\omega,x,y)=d(T^{k}%
\omega,x,y).
\]

Let us verify the assumptions in (\textbf{A2}). To check (a) observe that the
set-valued mapping $\omega\mapsto Y(\omega)$ is $\mathcal{F}_{0}$-measurable
because its graph is the intersection of the $\mathcal{F}_{0}\times
\mathcal{V}$-measurable sets $\{(\omega,x):\langle\phi(\omega),x\rangle=1\}$
and $\{(\omega,x):x\in K^{\circ}(\omega)\}$ (as regards the second set, see
Proposition A.3 in the Appendix).

To verify (b) consider a real number $r$ and the set
\[
Q=\{(\omega,x,y):x,y\in Y(\omega),\ d(\omega,x,y)>r\}.
\]
We have to show that $Q\in$ $\mathcal{F}_{0}\times\mathcal{V}\times
\mathcal{V}$. To this end observe that $d(\omega,x,y)>r$ if and only if%
\[
\inf_{j}\{\beta_{j}:x\leq_{K(\omega)}\beta_{j}y\}>e^{r}\sup_{j}\{\alpha
_{j}:\alpha_{i}y\leq_{K(\omega)}x\},
\]
where $\alpha_{j}>0$ and $\beta_{j}>0$ are rational numbers. By combining this
observation with the fact that $\{(\omega,x,y):x\leq_{K(\omega)}%
y\}\in\mathcal{F}_{0}\times\mathcal{V}\times\mathcal{V}$ (following from the
$\mathcal{F}_{0}$-measurability of $\omega\mapsto K(\omega)$) we obtain (b).

As we have noticed in Section 2, $(Y(\omega),d(\omega))$ is a complete
separable metric space and the topology generated by the metric $d(\omega)$ on
$Y(\omega)$ coincides with the Euclidean topology on $Y(\omega)$. From the
fact that $D(\omega,x)$ is completely monotone, it follows that the map
$f(\omega,x)$ transforms $Y(\omega)$ into $Y_{1}(\omega)$. Furthermore,
$f(\omega,\cdot)$ is continuous in the Euclidean topology and hence with
respect to the metric $d(\omega)$ on $Y(\omega)$ and the metric $d_{1}%
(\omega)$ on $Y_{1}(\omega)$. Therefore, conditions (c) and (d) hold.

Consider the $\mathcal{F}_{0}$-measurable sets $\Gamma_{1}\subseteq\Gamma
_{2}\subseteq...\subseteq\Omega$ constructed in Lemma 1. Let $\Omega
_{m}=\Gamma_{m+1}$,\ ($m=0,1,...$). We will show that the sets $\Omega
_{0}\subseteq\Omega_{1}\subseteq...\subseteq\Omega$ possess the properties
listed in (\textbf{A3}). Consider the mappings $f_{m}(\omega,x)$ and
$f^{(m)}(\omega,x)$ defined by (\ref{f1}) and (\ref{f2}), respectively. By
virtue of (\ref{f1}) and (\ref{f}), we get
\[
f_{m}(\omega,x)=f(T^{m-1}\omega,x)=\frac{D(T^{m-1}\omega,x)}{\langle\phi
_{m}(\omega),D(T^{m-1}\omega,x)\rangle},x\in X_{m-1}(\omega).
\]
Let us prove by induction with respect to $m=0,1,2,...$ the following formula
for every $x\in X_{-m-1}(\omega)$:%

\begin{equation}
f^{(m)}(\omega,x)=f_{0}(\omega)f_{-1}(\omega)...f_{-m}(\omega)x=\frac
{C(m+1,T^{-m-1}\omega)x}{\langle\phi(\omega),C(m+1,T^{-m-1}\omega)x\rangle}.
\label{fm-C}%
\end{equation}
If $m=0$, then
\[
f_{0}(\omega,x)=f(T^{-1}\omega,x)=\frac{D(T^{-1}\omega,x)}{\langle\phi
(\omega),D(T^{-1}\omega,x)\rangle}=\frac{C(1,T^{-1}\omega)x}{\langle
\phi(\omega),C(1,T^{-1}\omega)x\rangle},\ x\in X_{-1}(\omega).
\]
Suppose equation (\ref{fm-C}) holds for $m-1$. To verify it for $m$ we take
$x\in X_{-m-1}(\omega)$, put
\[
z=\frac{D(T^{-m-1}\omega)x}{\langle\phi_{-m}(\omega),D(T^{-m-1}\omega
)x\rangle}%
\]
and write%
\[
\frac{C(m+1,T^{-m-1}\omega)x}{\langle\phi(\omega),C(m+1,T^{-m-1}%
\omega)x\rangle}=\frac{C(m,T^{-m}\omega)D(T^{-m-1}\omega)x}{\langle\phi
(\omega),C(m,T^{-m}\omega)D(T^{-m-1}\omega)x\rangle}%
\]%
\[
=\frac{C(m,T^{-m}\omega)z}{\langle\phi(\omega),C(m,T^{-m}\omega)z\rangle
}=f_{0}(\omega)f_{-1}(\omega)...f_{-m+1}(\omega)z
\]%
\[
=f_{0}(\omega)f_{-1}(\omega)...f_{-m+1}(\omega)\frac{D(T^{-m-1}\omega
)x}{\langle\phi_{-m}(\omega),D(T^{-m-1}\omega)x\rangle}%
\]%
\[
=f_{0}(\omega)f_{-1}(\omega)...f_{-m+1}(\omega)f_{-m}(\omega)x.
\]
In this chain of equalities, the first one follows from (\ref{C-D}), the
second from the definition of $z$ and homogeneity of the mappings under
consideration, the third from the assumption of induction, the fourth from the
definition of $z$ and the last from the definition of $f_{-m}$.

Let $m$ be any nonnegative integer and let $\omega\in\Omega_{m}=\Gamma_{m+1}$.
Then, according to Lemma 1, the mapping $C(m+1,T^{-m-1}\omega)$ is strictly
monotone. Consequently (see (\ref{fm-C})), $f^{(m)}(\omega,x)\in K^{\circ
}(\omega)$ for all $x\in X_{-m-1}(\omega)$, and so $f^{(m)}(\omega
,X_{-m-1}(\omega))$ is a subset in $Y(\omega)$. Moreover, $f^{(m)}%
(\omega,X_{-m-1}(\omega))$ is a compact set in $Y(\omega)$, as a continuous
image of the set $X_{-m-1}(\omega)$ which is compact in the Euclidean
topology, and hence $f^{(m)}(\omega,X_{-m-1}(\omega))$ is compact with respect
to the metric $d(\omega)$. By virtue of Theorem 2, $f^{(m)}(\omega
,x):Y_{-m-1}(\omega)\rightarrow Y(\omega)$ is strictly non-expansive mapping
with respect to the metric $d_{-m-1}(\omega)$ on $Y_{-m-1}(\omega)$ and
$d(\omega)$ on $Y(\omega)$. Consequently, condition (\textbf{A3}) holds. Thus,
all the conditions sufficient for the validity of Theorem 3 are verified.

By virtue of assertion (i) of Theorem 3, there exists an $\mathcal{F}_{0}%
$-measurable mapping $\xi(\omega)\in Y(\omega)$ for which the equation
$\xi(T\omega)=f(\omega,\xi(\omega))\ $(a.s.) holds, i.e. we have
\[
\xi(T\omega)=\frac{D(\omega,\xi(\omega))}{\langle\phi_{1}(\omega),D(\omega
,\xi(\omega))\rangle}\ \text{(a.s.)}.
\]
Let $x(\omega)=\xi(\omega)$ and $\alpha(\omega)=\langle\phi_{1}(\omega
),D(\omega,\xi(\omega))\rangle$. Then $\alpha(\omega)x(T\omega)=D(\omega
)x(\omega)$ (a.s.) and since $x(\omega)\in Y(\omega)$, we have $x(\omega)\in
K^{\circ}(\omega)$, $\langle\phi(\omega),x(\omega)\rangle=1$ and
$\alpha(\omega)>0$. Furthermore, $x(\omega)$ is $\mathcal{F}_{0}$-measurable
and $\alpha(\omega)$ is $\mathcal{F}_{1}$-measurable, which proves assertion (a).

To prove (b), take any $(\alpha^{\prime}(\omega),x^{\prime}(\omega))$, where
$\alpha^{\prime}(\omega)\geq0$, $x^{\prime}(\omega)\in K(\omega)$,
$\langle\phi(\omega),x^{\prime}(\omega))\rangle=1$, satisfying $\alpha
^{\prime}(\omega)x^{\prime}(T\omega)=D(\omega)x^{\prime}(\omega)$ (a.s.).
Then
\[
\alpha^{\prime}(\omega)=\alpha^{\prime}(\omega)\langle\phi_{1}(\omega
),x^{\prime}(T\omega)\rangle=\langle\phi_{1}(\omega),D(\omega)x^{\prime
}(\omega)\rangle>0,
\]
and so%
\[
x^{\prime}(T\omega)=\frac{D(\omega)x^{\prime}(\omega)}{\alpha^{\prime}%
(\omega)}=\frac{D(\omega)x^{\prime}(\omega)}{\langle\phi_{1}(\omega
),D(\omega)x^{\prime}(\omega)\rangle}=f(\omega,x^{\prime}(\omega))\text{
(a.s.)}.
\]
By virtue of assertion (ii) of Theorem 3, we have $x^{\prime}(\omega
)=x(\omega)$ (a.s.), and consequently, $\alpha^{\prime}(\omega)=\alpha
(\omega)$ (a.s.), which proves (b).

To prove (c) we observe that $x:=a/\langle\phi_{-m-1}(\omega),a\rangle\in
X_{-m-1}(\omega)$ for any $0\neq a\in K_{-m-1}(\omega)$, and by virtue of
(\ref{fm-C}), the following equations hold:
\[
\frac{C(m+1,T^{-m-1}\omega)a}{\langle\phi(\omega),C(m+1,T^{-m-1}%
\omega)a\rangle}=\frac{C(m+1,T^{-m-1}\omega)x}{\langle\phi(\omega
),C(m+1,T^{-m-1}\omega)x\rangle}=f^{(m)}(\omega,x).
\]
By using (\ref{x-converg}), we get%
\[
\lim_{m(\omega)\leq m\rightarrow\infty}\sup_{0\neq a\in K_{-m-1}(\omega
)}d(\omega,\xi(\omega),\frac{C(m+1,T^{-m-1}\omega)a}{\langle\phi
(\omega),C(m+1,T^{-m-1}\omega)a\rangle})=
\]%
\[
\lim_{m(\omega)\leq m\rightarrow\infty}\sup_{x\in X_{-m-1}(\omega)}%
d(\omega,\xi(\omega),f^{(m)}(\omega,x))=0\ \text{(a.s.)}.
\]
In view of (\ref{inequality}), we can replace here the Hilbert-Birkhoff metric
$d(\omega,x,y)$ by the norm $\Vert x-y\Vert$, which yields (\ref{lim}).

The proof is complete.

\section{Appendix}

Let $(\Omega,\mathcal{F},P)$ be a probability space such that the $\sigma
$-algebra $\mathcal{F}$ is complete with respect to the measure $P$. Let $V$
be a complete separable metric space and $\mathcal{V}$ its Borel $\sigma
$-algebra. Let $\omega\mapsto A(\omega)$ be a set-valued mapping assigning a
non-empty set $A(\omega)\subseteq V$ to each $\omega\in\Omega$. A measurable
mapping $\sigma:(\Omega,\mathcal{F})\rightarrow(V,\mathcal{V})$\emph{\ }is
said to be a measurable selection of $A(\omega)$ if $\sigma(\omega)\in
A(\omega)$ for every $\omega$.

\textbf{Proposition A.1. }(\cite{Castaing}, Theorem III.22.)\textit{ Let
}$\omega\mapsto A(\omega)$\textit{ be an $\mathcal{F}$-measurable set-valued
mapping assigning a non-empty set }$A(\omega)\subseteq V$\textit{ to each
}$\omega\in\Omega$.\textit{ Then there exists a sequence of $\mathcal{F}%
$-measurable selections }$\sigma_{1}(\omega),\sigma_{2}(\omega),...$\textit{
such that for each }$\omega$\textit{ the sequence }$\{\sigma_{n}(\omega
)\}$\textit{ is dense in }$A(\omega)$\textit{.}

For a set $U\subseteq V$ define $A^{-1}(U)=\{\omega\in\Omega:A(\omega)\cap
U\neq\emptyset\}$ and denote by $\delta(x,U)$ the distance from $x$ to $U$.

\textbf{Proposition A.2. }(\cite{Castaing}, Theorem III.30.)\textit{ Let
}$\omega\mapsto A(\omega)$\textit{ be a set-valued mapping assigning a
non-empty closed set }$A(\omega)\subseteq V$\textit{ to each }$\omega\in
\Omega$\textit{. The following properties are equivalent:}

\textit{(a) The relation }$A^{-1}(U)\in\mathcal{F}$\textit{ holds either for
all open sets }$U$\textit{, or for all closed sets }$U$\textit{, or for all
Borel sets }$U$\textit{ in }$V$\textit{.}

\textit{(b) }$\delta(x,A(\cdot))$\textit{ is measurable for every }$x\in
V$\textit{.}

\textit{(c) The set-valued mapping }$A(\omega)$\textit{ admits a sequence of
measurable selections }$\sigma_{1}(\omega),\sigma_{2}(\omega),...$\textit{
such that for each }$\omega$\textit{ the sequence }$\{\sigma_{n}(\omega
)\}$\textit{ is dense in }$A(\omega)$\textit{.}

\textit{(d) The graph }$\{(\omega,a)\in\Omega\times V:a\in A(\omega)\}$
\textit{of the set-valued mapping }$A(\omega)$\textit{ belongs to
}$\mathcal{F\times V}$.

\textbf{Proposition A.3.}\textit{ Let }$\omega\mapsto K(\omega)$\textit{ be an
}$\mathcal{F}$\textit{-measurable set-valued mapping, where for each }%
$\omega\in\Omega$\textit{, the set }$K(\omega)$\textit{ is a solid cone in a
finite-dimensional vector space }$V$\textit{. Then the following assertions
are valid:}

\textit{(i) The set-valued mapping }$\omega\mapsto K^{\ast}(\omega)$\textit{
is }$\mathcal{F}$\textit{-measurable.}

\textit{(ii) The set-valued mappings }$\omega\mapsto V\setminus K^{\circ
}(\omega)$ and $\omega\mapsto K^{\circ}(\omega)$\textit{ are }$\mathcal{F}%
$\textit{-measurable.}

\textit{Proof. }Let us denote by $\mathcal{V}^{\ast}$ the Borel $\sigma
$-algebra on $V^{\ast}$. We wish to verify that the graph $\{(\omega,\phi
)\in\Omega\times V^{\ast}:\phi\in K^{\ast}(\omega)\}$ of the set-valued
mapping $K^{\ast}(\omega)$ belongs to $\mathcal{F}\times\mathcal{V}^{\ast}$.
By definition, $\phi\in K^{\ast}(\omega)$ if and only if $\langle\phi
,x\rangle\geq0\ \text{for every }x\in K(\omega)$. By assertion (c) of
Proposition (A.2), there exists a countable dense set $\{x_{n}(\omega
)\}_{n=1,2,...}$ of measurable selections of $K(\omega)$. Then $\langle
\phi,x\rangle\geq0\ $for every $x\in K(\omega)$ if and only if $\inf
_{n}\langle\phi,x_{n}(\omega)\rangle\geq0$, and so%
\[
\{(\omega,\phi)\in\Omega\times V^{\ast}:\phi\in K^{\ast}(\omega)\}=\{(\omega
,\phi):\inf_{n}\langle\phi,x_{n}(\omega)\rangle\geq0\}\in\mathcal{F}%
\times\mathcal{V}^{\ast}.
\]

Let $H(\omega)=V\setminus K^{\circ}(\omega)$. Let us show that the set-valued
mapping $\omega\mapsto H(\omega)$ is $\mathcal{F}$-measurable. By virtue of
assertion (b) of Proposition A.2, it is sufficient to verify that $\{\omega
\in\Omega:H(\omega)\cap B=\emptyset\}\in\mathcal{F}$ for every closed ball $B$
in $V$. This is so because%
\[
\left\{  \omega:H(\omega)\cap B=\emptyset\right\}  =\left\{  \omega:B\subseteq
K^{\circ}(\omega)\right\}  =\bigcup\limits_{k=1}^{\infty}\left\{
\omega:B_{[1\backslash k]}\subseteq K(\omega)\right\}
\]%
\[
=\bigcup\limits_{k=1}^{\infty}\bigcup\limits_{m=1}^{\infty}\left\{
\omega:b_{m}^{k}\in K(\omega)\right\}  =\bigcup\limits_{k=1}^{\infty}%
\bigcup\limits_{m=1}^{\infty}\{\omega:\delta(b_{m},K(\omega))=0\},
\]
where $B_{[1\backslash k]}$ is the closed $1/k$-neighborhood of $B$ and
$\{b_{m}^{k}\}_{m=1,2,...}$ is a dense subset in $B_{[1\backslash k]}$.
Therefore, the set $\{\omega:H(\omega)\cap B=\emptyset\}$ is a union of a
countable family of $\mathcal{F}$-measurable sets and thus belongs to
$\mathcal{F}$.

Now let us show that the set-valued mapping $\omega\mapsto K^{\circ}(\omega)$
is $\mathcal{F}$-measurable. Clearly, $a\in K^{\circ}(\omega)$ if and only if
$\delta(a,H(\omega))>0$. Consequently, we have%

\[
\left\{  (\omega,a)\in\Omega\times V:a\in K^{\circ}(\omega)\right\}  =\left\{
(\omega,a):\delta(a,H(\omega))>0\right\}  .
\]
Since $\omega\mapsto H(\omega)$ is $\mathcal{F}$-measurable, $\delta
(a,H(\omega))$ is $\mathcal{F}$-measurable by virtue of part (b) of
Proposition (A.2). Furthermore, $\delta(a,H(\omega))$ is continuos in $a$, and
so $\delta(a,H(\omega))$ is $\mathcal{F}\times\mathcal{V}$-measurable. The
proof is complete.



\begin{thebibliography}{99}                                                                                               %


\bibitem {Arnold}{Arnold L.}, \textit{Random dynamical systems}. {Springer},
Berlin (1998).

\bibitem {Arnold1994}Arnold L., Gundlach V.M. and Demetrius L., Evolutionary
formalism for products of positive random matrices. \textit{Ann. Appl. Prob.}
\textbf{4}, 859-901 (1994).

\bibitem {Arnlod-Evst-Gund}Arnold L, Evstigneev I.V. and Gundlach V.M.,
Convex-valued random dynamical systems: A variational principle for
equilibrium states. \textit{Random Oper. Stoc. Eqs.} \textbf{7}, 23-38 (1999).

\bibitem {Birkhoff1957}{Birkhoff G.}, {Extensions of Jentzsch's theorem}.
\textit{Trans. Amer. Math. Soc.} \textbf{84,} 219-227 (1957).

\bibitem {Castaing}{Castaing C. and Valadier M., \textit{Convex analysis and
measurable Multifunctions}. Lecture Notes in Mathematics, No 580,
Springer-Verlag, Berlin, Heidelberg, New York (1977).}

\bibitem {Dellacherie}{Dellacherie C. and Meyer P.-A.,} \textit{Probabilities
and Potential}. North Holland, Amsterdam (1978).

\bibitem {Dempster-Evstigneev}Dempster M.A.H., Evstigneev I.V. and
Schenk-Hopp\'{e} K.R., Exponential growth of fixed-mix strategies in
stationary asset markets. \textit{Finance Stoc.} \textbf{7}, 263-276 (2003).

\bibitem {DET}Dempster M.A.H., Evstigneev I.V. and Taksar M.I., Asset pricing
and hedging in financial markets with transaction costs: An approach based on
the von Neumann-Gale model. \textit{Annals of Finance} \textbf{2}, 327-355 (2006).

\bibitem {Dobrushin}Dobrushin R.L., Central limit theorem for nonstationary
Markov chains, I and II. \textit{Theory Prob. Appl.} \textbf{1}, 65-80,
329-383 (1956).

\bibitem {Eisenack}{Eisenack G. and Fenske C.}, \textit{Fixpunkttheorie}.
{BI-Wissenschaftsverlag,} Mannheim (1978).

\bibitem {Evstigneev1974}Evstigneev I.V., Positive matrix-valued cocycles over
dynamical systems. \textit{Uspekhi Matem. Nauk (Russian Mathematical Surveys)}
\textbf{29}, 219-220 (1974) (in Russian).

\bibitem {ESH-growth-model}{Evstigneev I.V. and }Schenk-Hopp\'{e} K.R., The
von Neumann-Gale growth model and its stochastic generalization. in: Dana
R.-A., Le Van C., Mitra T. and Nishimura K., (eds.), \textit{Handbook on
Optimal Growth}, Springer, New York \textbf{1, }337-383 (2006).

\bibitem {Evstigneev2007}{Evstigneev I.V. and Pirogov S.A.}, {A stochastic
contraction principle}. \textit{Random Oper. Stoc. Eqs.} \textbf{15}, 155-162 (2007).

\bibitem {ESH}Evstigneev I.V. and Schenk-Hopp\'{e} K.R., Stochastic equilibria
in von Neumann-Gale dynamical systems. \textit{Trans. Am. Math. Soc.}
\textbf{360}, 3345-3364 (2008).

\bibitem {Evstigneev2010}{Evstigneev I.V. and Pirogov S.A.}, {Stochastic
nonlinear Perron-Frobenius theorem}. \textit{Positivity} \textbf{14}(1), 43-57 (2010).

\bibitem {Frobenius1}Frobenius G., {\"{U}ber Matrizen aus positiven Elementen.
\textit{S.-B. Preuss. Akad. Wiss (Berlin)}, 471-6 (1908).}

\bibitem {Frobenius2}Frobenius G., {\"{U}ber Matrizen aus positiven Elementen,
II. \textit{S.-B. Preuss. Akad. Wiss (Berlin)}, 514-18 (1909).}

\bibitem {Frobenius3}Frobenius G., {\"{U}ber Matrizen aus nicht negativen
Elementen. \textit{S.-B. Preuss. Akad. Wiss (Berlin)}, 456-77 (1912).}

\bibitem {Fujimoto}Fujimoto T., A generalization of the Perron-Frobenius
theorem to nonlinear positive operators in a Banach space. \textit{The Kagawa
Univ. Eco.} \textbf{59}(3), 141-150 (1986).

\bibitem {Gaubert}Gaubert S. and Gunawardena J., The Perron-Frobenius theorem
for homogeneous, monotone functions. \textit{Trans. Am. Math. Soc.}
\textbf{356}, 4931-4950 (2004).

\bibitem {Hilbert}{Hilbert D.}, {\"{U}ber die gerade Linie als k\"{u}rzeste
Verbindung zweier Punkte}. \textit{Math. Ann.} \textbf{46}, 91-96 (1895).

\bibitem {Kifer1996}Kifer Yu., Perron-Frobenius theorem, large deviations, and
random perturbations in random environments. \textit{Math. Z.} \textbf{222},
677-698 (1996).

\bibitem {Kifer2008}Kifer Yu., Thermodynamic formalism for random
transformations revisited. \textit{Stoc. Dyn.} \textbf{8}, 77-102 (2008).

\bibitem {Kohlberg}{Kohlberg E.}, {The Perron-Frobenius theorem without
additivity}. \textit{J. Math. Econ.} \textbf{10}, 299-303 (1982).

\bibitem {Kolmogorov}Kolmogorov A.N. and Fomin S.V., \textit{Elements of the
Theory of Functions and Functional Analysis}. Graylock, Rochester, NY (1957).

\bibitem {Lemmens2}{Lemmens B. and Nussbaum R.}, \textit{Nonlinear
Perron-Frobenius theory}. {Cambridge Tracts in Mathematics 189, Cambridge
Univ. Press, Cambridge} (2012).

\bibitem {Morishima}Morishima M., \textit{Equilibrium, stability and growth}.
Clarendon Press, Oxford, (1964).

\bibitem {Nussbaum1}Nussbaum R.D., Hilbert's projective metric and iterated
nonlinear maps. \textit{Mem. Amer. Math. Soc.} \textbf{75}(391), (1988).

\bibitem {Nussbaum2}Nussbaum R.D., Iterated nonlinear maps and Hilbert's
projective metric II. \textit{Mem. Amer. Math. Soc.} \textbf{79}(401), (1989).

\bibitem {Nikaido}Nikaido H., \textit{Convex structures and economic theory}.
Academic Press, New York, (1968).

\bibitem {Oshime1}Oshime Y., Perron-Frobenius problem for weakly sublinear
maps in a Euclidean positive orthant. \textit{Japan J. Indust. Appl. Math.}
\textbf{9}(2), 313-350 (1992).

\bibitem {Oshime2}Oshime Y., An extension of Morishima's nonlinear
Perron-Frobenius theorem. \textit{J. Math. Kyoto Univ.} \textbf{23}, 803-830 (1983).

\bibitem {Oshime3}Oshime. Y, Non-linear Perron-Frobenius problem for weakly
contractive transformations. \textit{Math. Japon.} \textbf{29}, 681-704 (1984).

\bibitem {Perron1}Perron O., Grundlagen f\"{u}r eine Theorie des Jacobischen
Kettenbruchalgorithmus. \textit{Math. Ann.} \textbf{64}, 1-76 (1907).

\bibitem {Perron2}Perron O., Zur Theorie der Matrices. \textit{Math. Ann.}
\textbf{64}, 248-263 (1907).

\bibitem {S-S}Solow R.M., Samuelson P.A., Balanced growth under constant
returns to scale. \textit{Econometrica} \textbf{21}, 412-424 (1953).
\end{thebibliography}
\end{document}